\documentclass[fleqn,11pt]{article}
\usepackage{hyperref}
\usepackage{amsmath, amssymb, multirow, graphicx,footnote,caption}
\usepackage{color}
\newtheorem{algorithm}{Algorithm}[section]

\begin{document}

\title{\bf{A New Non-Linear Conjugate Gradient Algorithm for Destructive Cure Rate Model and a Simulation Study: Illustration with Negative Binomial Competing Risks}}
\author{\bf{Suvra Pal${}^{1}$\footnote{Corresponding author. E-mail address: suvra.pal@uta.edu Tel.: 817-272-7163.\newline \indent. } and Souvik Roy${}^{1}$}\\\\
${}^{1}$Department of Mathematics, University of Texas at Arlington, TX, 76019, USA.}
\date{}

\maketitle

\begin{abstract}

\noindent In this paper, we propose a new estimation methodology based on a projected non-linear conjugate gradient (PNCG) algorithm with an efficient line search technique. We develop a general PNCG algorithm for a survival model incorporating a proportion cure under a competing risks setup, where the initial number of competing risks are exposed to elimination after an initial treatment (known as destruction). In the literature, expectation maximization (EM) algorithm has been widely used for such a model to estimate the model parameters. Through an extensive Monte Carlo simulation study, we compare the performance of our proposed PNCG with that of the EM algorithm and show the advantages of our proposed method. Through simulation, we also show the advantages of our proposed methodology over other optimization algorithms (including other conjugate gradient type methods) readily available as R software packages. To show these we assume the initial number of competing risks to follow a negative binomial distribution although our general algorithm allows one to work with any competing risks distribution. Finally, we apply our proposed algorithm to analyze a well-known melanoma data.


\end{abstract}


\section{Introduction}

Due to the advancements in the treatment of certain types of disease, including cancer and heart disease, we see a significant number of patients to respond favorably to the treatment and not show recurrence until the end of a long follow-up time. In literature, these patients are called recurrence-free survivors. It may be possible that some of these recurrence-free survivors will not show recurrence for a sufficiently long period after the follow-up time since they may reach a stage where the disease is undetectable as well as harmless. These patients, among the recurrence-free survivors, are called long-term survivors or ``cured''. It is to be noted that the estimation of this long-term survivor rate or cure rate cannot be readily obtained from a given survival data since we are not in a position to identify which of the recurrence-free survivors can be considered as long-term survivors. This issue arises because a patient who is susceptible to disease recurrence soon after the follow-up time may also show no recurrence and survive until the end of the follow-up time. However, the estimation of a treatment-specific cure rate is crucial to observe the trend in the survival of patients suffering from a particular disease. Furthermore, it is an important measure to judge the efficacy of a treatment and its adoption in practice, as opposed to the standard treatment.

The literature on cure rate models is vast and the topic itself is one of the most emerging areas of modern research. The early work on cure rate model dates back to the work of Boag (1949) followed by Berkson and Gage (1952), which is known in the literature as the mixture cure rate model. According to the mixture cure rate model, the overall survival function (also called the population survival function) of the time-to-event variable $Y$ can be split into two parts, one corresponding to the cured group and the other corresponding to the susceptible group. Such a survival function is given by
\begin{equation}
S_{pop}(y) = p_0 + (1-p_0)S_s(y),
\label{mix}
\end{equation}
where $p_0$ is the proportion of cured subjects (cure rate) and $S_s(y)$ is a proper survival function for the susceptible subjects; see Sy and Taylor (2000) and Kannan et al. (2010). Note that $S_{pop}(y)$ is not a proper survival function as $\lim_{y\rightarrow\infty}S_{pop}(y) = p_0$ $(\neq 0)$. One major drawback of the mixture cure rate model in \eqref{mix} is that 
it does not incorporate a scenario where several risk factors may compete to produce the event of interest; known as the competing risks scenario. To circumvent this problem, Chen et al. (1999) proposed the promotion time cure rate model by considering a competing risks scenario and assuming the latent number of risk factors to follow a Poisson distribution. The corresponding population survival function is given by 
\begin{equation}
S_{pop}(y) = e^{-\eta(1-S(y))},
\label{prom}
\end{equation}
where $\eta$ is the mean number of risk factors and $S(y)$ is the common survival function of the progression times, defined as the time taken by each risk factor to produce the event. Note that in this case the cure rate is given by $e^{-\eta}$. Rodrigues et al. (2009) unified the mixture and promotion time cure rate models by proposing the Conway-Maxwell Poisson (COM-Poisson) cure rate model, which assumes the number of risk factors to follow a COM-Poisson distribution that can handle both over-dispersion and under-dispersion. For this unified cure rate model, the population survival function is given by
\begin{equation}
S_{pop}(y)=\frac{Z(\eta S(y),\phi)}{Z(\eta,\phi)},
\label{com}
\end{equation}
where $Z(a,\phi) = \sum_{j=0}^\infty{\frac{a^{j}}{(j!)^{\phi}}}$ and $\phi$ is the dispersion parameter of the COM-Poisson distribution. In \eqref{com}, $\eta$ is related to the mean number of risk factors and $S(y)$ is as defined before for the promotion time model. In this case, the cure rate is given by $\frac{1}{Z(\eta,\phi)}$. Note that if $\phi\rightarrow\infty$ in \eqref{com}, the COM-Poisson model reduces to the mixture model in \eqref{mix} with $p_0 = \frac{1}{1+\eta}$, whereas, if $\phi=1$ in \eqref{com}, the COM-Poisson model reduces to the promotion time model in \eqref{prom}; see Cancho et al. (2011), Balakrishnan and Pal (2015), Balakrishnan and Pal (2016), and Balakrishnan and Feng (2018) for some recent works on cure rate model using COM-Poisson distribution. To develop the associated inferential procedures, several approaches have been proposed in the literature. In this regard, one may refer to parametric approaches (Farewell, 1986; deFreitas and Rodrigues, 2013; Balakrishnan and Pal, 2013); semi-parametric approaches (Kuk and Chen, 1992; Li and Taylor, 2002; Balakrishnan et al., 2017); and non-parametric approaches (Maller and Zhou, 1996; Balakrishnan et al., 2016). Rodrigues et al. (2011) first brought in a practical and interesting interpretation of the biological mechanism of the occurrence of an event of interest; see Cooner et al. (2007). They proposed a flexible cure rate model, known as the destructive cure rate model, by considering the possible elimination (or destruction) of risk factors after an initial treatment. Since then, several papers have been published in the context of destructive model and interested readers may refer to Cancho et al. (2013), Pal and Balakrishnan (2016), Pal et al. (2018), and Gallardo et al. (2016), among others.

In this paper, we consider a competing risks scenario and accommodate the elimination of risk factors after an initial treatment. The corresponding model is termed as the destructive cure rate model and such a model was first proposed by Rodrigues et al. (2011), where the authors assumed a weighted Poisson distribution for the initial number of competing risks and carried out a direct maximization of the observed likelihood function for parameter estimation. Since then, few other authors studied the destructive cure rate model and interested readers may refer to Pal and Balakrishnan (2016), Gallardo et al. (2016), Pal and Balakrishnan (2017), and Pal et al. (2018), among others. In particular, Pal and Balakrishnan (2016), Pal and Balakrishnan (2017), and Pal et al. (2018) studied the destructive cure rate model by assuming different distributions for the initial risk factors. For the maximum likelihood estimation (MLE) of the model parameters, the authors developed the expectation maximization (EM) algorithm; see McLachlan and Krishnan (2008). However, they noted that the likelihood surface was flat with respect to certain parameters of the model and as such simultaneous maximization of all model parameters was an issue. To circumvent this problem, the authors proposed profile likelihood techniques within the EM algorithm. Although the proposed approach performed satisfactorily, few drawbacks were noted. For instance, the root mean square error (RMSE) of the regression parameters associated with the cure rate turned out to be high. Such inaccuracy may lead to imprecise inference on the overall population survival. Furthermore, the suggested profile likelihood techniques require the EM algorithm to be run several times which is computationally expensive. To evade these issues with the developed EM algorithm, we propose a new estimation procedure based on a projected non-linear conjugate gradient (PNCG) algorithm with an efficient line search (Hager and Zhang, 2005) that (i) allows simultaneous maximization of all model parameters; (ii) results in more precise estimates, specifically for the parameters associated with the cure rate; and (iii) is computationally less expensive. We also compare our proposed PNCG algorithm with the EM algorithm as well as with other optimization algorithms, including available conjugate type methods, readily available as R software packages and show the main advantages of using our methodology. To the best of our knowledge, we are the first one to propose such an algorithm in the context of cure rate models. 

The rest of the paper is organized as follows. In Section 2, we give a brief overview of the destructive cure rate model. In Section 3, we describe our proposed PNCG algorithm in detail. In Section 4, we present the results of a detailed simulation study where we assume the initial competing risks to follow a negative binomial distribution. We first compare the performance of our PNCG algorithm with the EM algorithm developed by Pal and Balakrishnan (2016). Then, we also compare the PNCG algorithm with other optimization techniques available as R packages. In Section 5, we apply the PNCG algorithm to analyze to a well-known melanoma data. Finally, in Section 6, we make some concluding remarks and discuss some future research in this direction.

\section{Destructive cure rate model}

We can define the destructive cure rate model as follows; see Rodrigues et al. (2011). First, we assume that there are $M$ latent risk factors competing to produce an event of interest (for instance, death due to cancer or recurrence of a disease). These risk factors being unobserved, we assume them to follow a discrete distribution with mass function $p_m = P[M=m]$. With the passage of time or after an initial treatment, we consider the possibility of elimination of risk factors, which we assume to occur according to a Binomial law. More specifically, we define the number of active risk factors after elimination, i.e., those risk factors that are still capable of producing the event, as
\begin{equation}
    D= \begin{cases}
              X_1+X_2+\ldots+X_M,&  M>0,\\
               0, &  M = 0,
           \end{cases}
\label{Des}
\end{equation}
where $X_j's$ are independent Bernoulli random variables with $P[X_j=1] = p,$ for $j=1,2,\cdots,M,$ with $p$ denoting the activation probability of each risk factor. Note that these assumptions are in line with the assumptions of Rodrigues et al. (2011) who first proposed the destructive cure rate model; see also Yang and Chen (1991). Using the distribution of $M$ and noting that the conditional distribution of $D$ given $M=m$ is Binomial $(m,p)$, the marginal distribution of $D$ can be obtained. Given $D=d,$ we now let $W_j$ $(j=1,2,\ldots,d)$ to denote the time taken by the $j-$th active risk factor to produce the event, also called the progression time. Once again, using the same assumptions as in Rodrigues et al. (2011), we let the progression times $W_j's$ to be independently distributed and distributed independently of $D$ with distribution function $F(\cdot)=1-S(\cdot),$ where $S(\cdot)$ is the corresponding survival function. Note that the individual progression times are not observed, however, we only observe the time taken by the first active risk factor to produce the event, which we term as the lifetime. Notationally, such a lifetime in a competing risks scenario is defined as
\begin{equation}
Y= \begin{cases}
\text{min}\{W_1,\ldots,W_D\}, & D > 0\\
\infty, & D = 0.
\end{cases}
\label{Y}
\end{equation}
The infinite lifetime corresponding to $D=0$ leads to a proportion, say $p_0$, of the population who are not susceptible to the occurrence of the event. We term this proportion as the ``cure rate" and its estimation is of great interest to us. As pointed out by Rodrigues et al. (2011), the destructive cure rate model is not identifiable as per Li et al. (2001). One way to circumvent this issue is to bring in the effect of prognostic factors (or covariates). For instance, we can relate the parameter $p$ to a set of covariates $\boldsymbol {x}$ using the logistic link function $p=\frac{\exp(\boldsymbol x^{'}\boldsymbol\beta_1)}{1+\exp(\boldsymbol x^{'}\boldsymbol\beta_1)}$ and a suitable parameter related to the distribution of $M$ to another set of covariates $\boldsymbol {z}$ using an appropriate link function $g(\boldsymbol z^{'}\boldsymbol\beta_2)$, where $\boldsymbol\beta_1$ and $\boldsymbol\beta_2$ represents the vectors of regression coefficients. Furthermore, we have to make sure that either $\boldsymbol\beta_1$ or $\boldsymbol\beta_2$ do not include the intercept term to retain identifiability. Note also that $\boldsymbol x$ and $\boldsymbol z$ cannot share common elements. The survival function, also called the long-term survival function, of the random variable $Y$ in \eqref{Y} is given by
\begin{equation*}
S_{pop}(y) = P[Y \geq y] = \sum_{d=0}^\infty P[D=d]\{S(y)\}^d.
\end{equation*}
The corresponding density function, called the long-term density function, can be easily obtained as $f_{pop}(y)=-S_{pop}^\prime (y)$. On assuming specific distributions for $M$, and hence for $D$, explicit expressions for $S_{pop}(y)$ and $f_{pop}(y)$ can be obtained.

\section{Estimation method: projected non-linear conjugate gradient algorithm with line search}
We consider a scenario where the lifetime data may not be completely observed and is thus subject to right censoring. We let $T_i$ and $C_i$ denote the actual failure time and censoring time, respectively, $i=1,2,\cdots,n,$ where $n$ denotes the sample size. The lifetime that we observed is then given by $Y_i=\min\{T_i,C_i\}$. We let $\delta_i$ denote the right censoring indicator, i.e., $\delta_i$ takes the value 1 if the lifetime is observed and 0 if it is right censored. The observed data can then be represented by $\boldsymbol O = \{(y_i,\delta_i,\boldsymbol x_i, \boldsymbol z_i), i=1,2,\cdots,n\}$. On assuming non-informative right censoring, we can define the observed data likelihood function as 
\begin{eqnarray}
L(\boldsymbol\theta)= \prod_{i=1}^n \{f_{pop}(y_i|\boldsymbol {x_i}, \boldsymbol {z_i})\}^{\delta_i}\{S_{pop}(y_i|\boldsymbol{x_i}, \boldsymbol {z_i})\}^{1-\delta_i}, \label{eq:likelihood}
\end{eqnarray}
where $\boldsymbol{\theta}$ is the vector of unknown parameters. We, now, desire to estimate the optimal parameter set, denoted by $\hat{\boldsymbol{\theta}}$, that maximizes the likelihood function given in (\ref{eq:likelihood}). Applying the natural logarithm on both sides of (\ref{eq:likelihood}), we obtain the log-likelihood function as follows
\begin{equation}\label{eq:loglikelihood}
l(\boldsymbol\theta)= \sum_{i=1}^n[ \delta_i \log\{f_{pop}(y_i|\boldsymbol {x_i}, \boldsymbol {z_i})\}+{(1-\delta_i)}\log\{S_{pop}(y_i|\boldsymbol{x_i}, \boldsymbol {z_i})\}].
\end{equation}
Then, the corresponding maximization problem to obtain the optimal parameter set is given by 
\begin{equation}\label{eq:maxproblem}
\hat{\boldsymbol\theta}={\arg\max}_{\boldsymbol\theta\in U}~ l(\boldsymbol\theta),\\
\end{equation}
where $U$ is the feasible set of constraints.

The function $l$ is non-linear with respect to $\boldsymbol\theta$ and, thus, gives rise to a non-linear maximization problem. To solve such a non-linear maximization problem (\ref{eq:maxproblem}), we use the non-linear conjugate gradient method (NCG) together with an efficient line search technique based on the formula of Hager and Zhang (2005). This method has been primarily used in the context of solving partial differential equation (PDE)-constrained optimal control problems arising in mathematical models of crowd motion, game theory and medical imaging (Roy et al., 2016; Roy et al., 2017; Roy et al., 2018; Adesokan et al., 2018), but not well explored in the context of cure rate models. The NCG scheme has significant advantages over traditional Newton-based schemes (for instance, the EM algorithm in Pal and Balakrishnan, 2016, where the maximization step was carried out using a one-step Newton Raphson method). For example, even though Newton's method potentially converges faster, it is highly expensive and time consuming to compute the Hessian in a Newton-based method, specially for optimization problems with large number of parameters and bigger sample size $n$. On contrary, in the NCG method, only the gradient is required to be evaluated, leading to a much faster convergence than the Newton-based schemes. This is also observed in the simulation studies presented in Section \ref{sec:simulation}.

To start the NCG scheme, we use an initial guess $\boldsymbol\theta_0$ for the parameter set. It has been observed that the performance of the NCG method remains the same for a large choice of initial guesses, which suggests the robustness of the method; see Roy et al., (2017) and Adesokan et al., (2018).  Due to the fact that the maximum rate of increase of a function is along the positive gradient direction, the initial guess is updated by moving in the search direction given by the gradient $\boldsymbol{d}_0=\boldsymbol{g}_0=\frac{\partial}{\partial\boldsymbol\theta} l(\boldsymbol\theta_0)$ of the function $l$.  In subsequent iterations, the search directions are recursively given by the formula 
\[
\boldsymbol{d}_{k+1} = \boldsymbol{g}_{k+1}+\xi_k \boldsymbol{d}_k,~ k=0,1,2,\cdots,
\]
where $\boldsymbol{g}_k = \frac{\partial}{\partial\boldsymbol\theta} l(\boldsymbol\theta_k)$ and $\xi_k$ is given by the formula of Hager-Zhang (Hager and Zhang, 2005) as follows
\begin{equation}\label{eq:step_HG}
\xi_k = \dfrac{1}{\boldsymbol{d}_k' \boldsymbol{w}_k}\left(\boldsymbol{w}_k-2\boldsymbol{d}_k\dfrac{\boldsymbol{w}_k' \boldsymbol{w}_k}{\boldsymbol{d}_k' \boldsymbol{w}_k}\right)' \boldsymbol{g}_{k+1}
\end{equation}
with $\boldsymbol{w}_k = \boldsymbol{g}_{k+1}-\boldsymbol{g}_k$. We update our parameter set $\boldsymbol\theta$ using a steepest ascent scheme given below
\begin{equation}\label{eq:update}
\boldsymbol\theta_{k+1} = \boldsymbol\theta_k+s_k \boldsymbol{d}_k,
\end{equation}
where $s_k>0$ is a steplength obtained through a line search algorithm. An accurate estimate of the steplength $s_k$ is crucial because a very large $s_k$ would result in deviating from the path of the maximizer, whereas a very small $s_k$ would lead to slow convergence of the NCG scheme. Thus, we deploy a line-search algorithm to obtain $s_k$ that uses the following Armijo condition (Annunziato and Borzi, 2013) of sufficient increase of $l$
\[
l(\boldsymbol\theta_k + s_k \boldsymbol{d}_k) \geq l(\boldsymbol\theta_k) + \lambda s_k ~\boldsymbol{d}_k' \boldsymbol{g}_k
\]
with $0 < \lambda < 1/2$; see Nocedal and Wright (1999). Such a method leads to an optimal steplength $s_k$ resulting in a fast and accurate optimal solution through the NCG. 

Traditional gradient-based schemes look for an optimum in a global space of parameters. But for the maximization problem given in (\ref{eq:maxproblem}), one needs to determine the maximizer that lies in the constraint set $U$. Thus, we need to ensure that the solution obtained at the end of each iterative step in the NCG scheme lies in the constraint set $U$. For this purpose, we use a projection step onto the constraint set $U$ that is applied to the parameter update step (\ref{eq:update}) through the following way
\[
\boldsymbol\theta_{k+1} = \mathbb{P}[\boldsymbol\theta_k+s_k \boldsymbol{d}_k],
\]
where 
$\mathbb{P}[\boldsymbol\theta]$ is the projected parameter set. The projection step ensures that the parameter values obtained in each iteration lies inside the constraint set $U$. We call our new scheme as the projected NCG (PNCG) scheme with line search. The scheme is terminated once the relative difference between two successive iterates is less than a specified tolerance level or the number of iterations exceed the maximum number of iterations. The algorithm of the PNCG scheme is summarized below.
\begin{algorithm}[PNCG Scheme with line search]\label{alg:NCG}\
\begin{enumerate}
\item Input: initial guess $\boldsymbol\theta_0$. Evaluate $\boldsymbol{d}_0 = \frac{\partial}{\partial\boldsymbol\theta}l(\boldsymbol\theta_0)$, index $k=0$, maximum $k=k_{max}$, tolerance = $tol$
\item While $(k<k_{max}),$ do
\item Set $\boldsymbol\theta_{k+1} = \mathbb{P}[\boldsymbol\theta_k + s_k\boldsymbol{d}_k]$, where $s_k$ is obtained using the line-search algorithm
\item Compute $\boldsymbol{g}_{k+1} = \frac{\partial}{\partial\boldsymbol\theta}l(\boldsymbol\theta_{k+1})$
\item Compute $\xi_k$ using (\ref{eq:step_HG})
\item Set $\boldsymbol{d}_{k+1}=\boldsymbol{g}_{k+1}+\xi_k \boldsymbol{d}_k$
\item If $\|\frac{\boldsymbol\theta_{k+1}-\boldsymbol\theta_k}{\boldsymbol\theta_k}\| < tol$, terminate
\item Set $k=k+1$
\item End while.
\end{enumerate}
\end{algorithm}
The convergence of the PNCG scheme, as described in Algorithm \ref{alg:NCG}, follows from Neittaanmaki and Tiba (1994, Lemma 1.6, p. 235).

\section{Simulation study with negative binomial competing risks}\label{sec:simulation}

We shall now assume the distribution of the initial number of competing risks to be negative binomial with the following mass function
\begin{equation}
p_m = P[M=m;\eta,\phi]=\frac{\Gamma(m+\frac{1}{\phi})}{\Gamma(\frac{1}{\phi})m!}\bigg(\frac{\phi\eta}{1+\phi\eta}\bigg)^m(1+\phi\eta)^{-\frac{1}{\phi}}, \ \  m=0,1,2,\ldots,
\label{nb}
\end{equation}
where $\eta>0$ and $\phi>0$. Using \eqref{nb} and noting that the conditional distribution of $D$ given $M=m$ is Binomial ($m$,$p$), Pal and Balakrishnan (2016) showed that the mass function of $D$ can be expressed as
\begin{eqnarray}
P[D=d;\eta,\phi,p] =
\frac{\Gamma\bigg(d+\frac{1}{\phi}\bigg)}{\Gamma\bigg(\frac{1}{\phi}\bigg)d!}\bigg(\frac{\phi\eta p}{1+\phi\eta p}\bigg)^d\bigg(\frac{1}{1+\phi\eta p}\bigg)^{\frac{1}{\phi}}, \ \ d=0,1,2,\cdots,
\label{D}
\end{eqnarray}
which is once again a negative binomial distribution. Rodrigues et al. (2011) showed that the overall survival function or the population survival function of the lifetime variable $Y$ in \eqref{Y} can be expressed as
\begin{eqnarray}
S_{pop}(y) = P[Y>y] = \{1+\phi\eta pF(y)\}^{-\frac{1}{\phi}}
\label{S_pop}
\end{eqnarray}
and the corresponding density function can be expressed as
\begin{equation}
f_{pop}(y)=\frac{\eta p}{1+\phi\eta pF(y)}S_{pop}(y)f(y),
\label{f_pop}
\end{equation}
where $f(\cdot)$ is the density function corresponding to $F(\cdot)$. Noting that cure rate is the long-term survival probability, we can easily obtain the cure rate as
\begin{equation}
p_0 = \bigg\{\frac{1}{1+\phi\eta p}\bigg\}^{\frac{1}{\phi}}.
\label{p0}
\end{equation}
In this regard, Pal and Balakrishnan (2016) showed that the overall cure rate can be decomposed into pre-destructive and post-destructive components. From \eqref{p0}, we note that the cure rate is a decreasing function of both $\eta$ and $p$. We will link the parameter $\eta$ using a log-linear link function $\eta=\exp(\boldsymbol z^{'}\boldsymbol\beta_2)$. Although we have assumed negative binomial distribution for competing risks, one can easily incorporate other competing risks distribution such as those considered in Pal and Balakrishnan (2017) and Pal et al. (2018), among others. As done in Pal and Balakrishnan (2016), we will also assume $F(\cdot)$ and $f(\cdot)$ to be the distribution function and density function, respectively, of a two-parameter Weibull distribution defined as
\begin{eqnarray}
F(y)&=&1-\exp\{-(\gamma_2 y)^{\frac{1}{\gamma_1}}\} \ \ \text{and} \nonumber \\ 
f(y) &=& \frac{1}{\gamma_1 y}(\gamma_2 y)^{\frac{1}{\gamma_1}}\{1-F(y)\} 
\label{wei}
\end{eqnarray}
for $y>0$, $\gamma_1>0,$ and $\gamma_2>0.$ Thus, we now have $\boldsymbol\theta=(\phi,\boldsymbol\beta_1^\prime,\boldsymbol\beta_2^\prime,\gamma_1,\gamma_2)^\prime$, $U = \lbrace \boldsymbol\theta : \phi>0, \boldsymbol{\beta_1},\boldsymbol{\beta_2},\gamma_1 > 0,\gamma_2 >0 \rbrace$, and
\[
\begin{aligned}
\mathbb{P}[\boldsymbol\theta] = &\lbrace \max(0,\phi),\boldsymbol\beta_1,\boldsymbol\beta_2, \max(0,\gamma_1),\max(0,\gamma_2)\rbrace.
\end{aligned}
\]

For the purpose of this simulation study, we mimic the real melanoma dataset that we analyze in the next section. For this dataset, two covariates of interest are tumor thickness (measured in mm) and ulceration status (presence of ulcer denoted by 1 and absence denoted by 0). A preliminary analysis of this data indicates that 44\% of patients had the presence of ulcer. For this group of patients, the mean and standard deviation of tumor thickness turned out to be 4.34 mm and 3.22 mm, respectively. On the other hand, for the group of patients without the presence of ulcer, the mean and standard deviation of tumor thickness turned out to be 1.81 mm and 2.19 mm, respectively. Moreover, as noted by Pal and Balakrishnan (2016), histograms of tumor thickness for two groups suggest that a Weibull distribution may be suitable for the group with presence of ulcer, whereas an exponential distribution may be suitable for the group without the presence of ulcer. Thus, to generate the ulceration status data, we first generated a Uniform (0,1) random variable, say, $U$. If $U \leq 0.44,$ we took ulceration status ($x$) to be 1 and then generated the tumor thickness ($z$) from a Weibull distribution. For the choice of the Weibull parameters here, we equated the theoretical mean and variance of the Weibull distribution to 4.34 and 10.37, respectively. On the other hand, if $U > 0.44,$ we took ulceration status ($x$) to be 0 and generated tumor thickness ($z$) from an exponential distribution. For a choice of this exponential parameter, we simply equated the exponential mean to 1.81. To make the model identifiable, we linked the parameter $\eta$ to ulceration status without the intercept term and the parameter $p$ to tumor thickness (including the intercept term). Thus, we have the following link functions: $\eta = \exp(\beta_2 z)$ and $p = \frac{\exp(\beta_{0}+\beta_{1}x)}{1+\exp(\beta_{0}+\beta_{1}x)}$. To decide on the true value of the regression parameter $\beta_2$, we note that the absence of ulcer would simply imply $\eta$ taking on the value 1. Intuitively, with the presence of ulcer, $\eta$ is expected to be higher since it is related to the mean number of active competing risks. Thus, we chose $\eta$ to be 3 in the presence of ulcer, which gave us the true value of $\beta_2$ as $\log(3) = 1.099.$ To decide on the true values of the regression parameters $\beta_0$ and $\beta_1$ corresponding to $p$, we chose a low and a high value of $p$ as 0.3 and 0.9, respectively.  Also, from the generated tumor thickness data, we identified the minimum ($x_{min}$) and the maximum ($x_{max}$) tumor thickness values. Then, we came up with the following two equations to solve for $\beta_0$ and $\beta_1$.
\begin{eqnarray*}
\frac{\exp(\beta_{0}+\beta_{1}x_{min})}{1+\exp(\beta_{0}+\beta_{1}x_{min})} &=& 0.3 \\
\frac{\exp(\beta_{0}+\beta_{1}x_{max})}{1+\exp(\beta_{0}+\beta_{1}x_{max})} &=& 0.9.
\end{eqnarray*}
Note that the true values of $\beta_0$ and $\beta_1$ depends on the generated tumor thickness data and as such cannot be kept fixed across the generated data sets; see Pal and Balakrishnan (2016). To incorporate random censoring, the censoring time $(C)$ distribution was chosen to be exponential with censoring rate 0.15. Next, to generate the lifetime data from the DNB model, we followed the following steps:
\begin{itemize}
\item[(i)] Generate the initial number of risk factors $M$ from a negative binomial distribution with mass function as in \eqref{nb} for a chosen value of $\phi.$
\item[(ii)] If $M = 0$ in (i), set the number of active risk factors $D = 0$.
\item[(iii)] If $M > 0$ in (i), generate $D$ from a Binomial distribution with the observed value of $M$ as the number of trials and success probability $p$.  
\item[(iv)] From (ii) and (iii), if $D = 0,$ set the observed lifetime $Y$ as the censoring time, i.e., $Y = C$. 
\item[(v)] From (ii) and (iii), if $D > 0,$ generate $D$ Weibull random variables $\{W_1,W_2,\cdots,W_D\}$ with density function as in \eqref{wei} for chosen values of $\gamma_1$ and $\gamma_2$. Then, take the observed lifetime as $Y = \min\{\min\{W_1,W_2,\cdots,W_D\},C\}$.
\item[(vi)] From (iv) and (v), if $Y = C,$ set $\delta = 0$, otherwise, set $\delta = 1$.
\end{itemize}
 As done in Pal and Balakrishnan (2016), we considered the following true choices: $(\gamma_1,\gamma_2)$ = (0.215,0.183) and (0.316,0.179); $n$ = 300 and 400; and $\phi$ = 0.5 and 0.75. We ran our simulations in R statistical software and all results were based on 500 Monte Carlo runs. To come up with a choice of initial values to start the PNCG algorithm, we first created an interval for each model parameter by taking 20\% deviation off its true value and then selected a value at random from the created interval, which was used as the parameter's initial value.
 
 To demonstrate the performance of our proposed PNCG algorithm with line search, we chose $k_{max} = 500,$ $\lambda = 0.1$, and $tol = 0.001.$ Note that we also tried the algorithm for other values of $\lambda$ within its range, and our findings were similar. First, we present the calculated bias and root mean square (RMSE) of the estimates obtained using PNCG algorithm and compare these with the bias and RMSE obtained using the EM algorithm developed by Pal and Balakrishnan (2016). In this regard, note that in the EM algorithm the maximization step was carried out using a one-step Newton Raphson method and the same stopping criterion as in the aforementioned PNCG algorithm was used. Furthermore, note that in the PNCG algorithm the maximization was done on the observed log-likelihood function, whereas in the EM algorithm the maximization was done on the $Q-$function (using the same notation as in Pal and Balakrishnan, 2016). 

\subsection{Comparison with EM algorithm}

In Tables \ref{table:phi0.5} and \ref{table:phi0.75}, we present the simulation study results, in terms of bias and RMSE when the true value of $\phi$ is 0.50 and 0.75, respectively. When we employed the PNCG algorithm, we did not face any issue with simultaneous maximization of all model parameters. This is a nice behavior of the PNCG algorithm and is unlike the EM algorithm developed by Pal and Balakrishnan (2016), where simultaneous maximization was not possible and the authors had to keep the parameters $\gamma_1$ and $\phi$ fixed for the EM procedure to work. This is a big advantage of the PNCG algorithm over the EM algorithm. It is clear that the PNCG algorithm provides estimates that are very close to the true parameter values. Both bias and RMSE are found to decrease with an increase in sample size, which is also a very satisfactory property. When compared to the estimates produced by EM algorithm, we first note that the bias in the estimates of the PNCG algorithm is lower than that of the EM algorithm. Note, in particular, the reduction in the bias of $\phi$ that the PNCG algorithm results when the true value of $\phi$ is 0.75. When compared to the EM algorithm, the PNCG algorithm also results in a significant reduction in the RMSE for the regression parameters $(\beta_{0}, \beta_{1},\beta_2)$ as well as for the parameter $\phi$. This is another big advantage of the PNCG algorithm over the EM algorithm, noting that the cure rate is a pure function of the regression parameters and $\phi$ only. For the lifetime parameters $\gamma_1$ and $\gamma_2,$ both algorithms result in similar bias and RMSE. From the above findings, it is very clear that the PNCG algorithm results in more accurate and precise estimates of the model parameters and is thus preferred over the EM algorithm.
\begin{table} [ht!]
\caption{Comparison of PNCG algorithm with EM algorithm in terms of bias and RMSE with true value of $\phi$ as 0.50}
\centering
\begin{tabular}{ l l l l l l l }\\ 
\hline                                 
$n$ & $(\gamma_1,\gamma_2)$ & Parameter &\multicolumn{2}{c}{Bias}  &\multicolumn{2}{c}{RMSE} \\  \cline{4-5}  \cline{6-7}
& & &PNCG & EM & PNCG & EM \\
\hline  

300 & (0.215,0.183)  & $\beta_{2}$                          & 0.021 & 0.035 & 0.229 &0.339  \\
 & &  $\beta_{0}$                                                   & -0.028  & -0.120 & 0.271 &0.462  \\
  & &  $\beta_{1}$                                                   & 0.015 & 0.115 &0.139  &0.396  \\
 & &  $\gamma_1$                                                    & -0.002 &-0.003  &0.019  &0.019  \\ 
 & &  $\gamma_2$                                                    & 0.001 & 0.000 &0.006  &0.007 \\
 & &  $\phi$                                                              & -0.012  & 0.012 &0.138  &0.370 \\[0.5ex]
 400 & (0.215,0.183)  & $\beta_{2}$                          & 0.019 & 0.024 &0.181  &0.309  \\
 & &  $\beta_{0}$                                                   & -0.018  & -0.070 & 0.201 & 0.379 \\
  & &  $\beta_{1}$                                                   & 0.013 & 0.053 &0.106  &0.253  \\
 & &  $\gamma_1$                                                    & -0.001 & -0.002 & 0.017 & 0.017 \\ 
 & &  $\gamma_2$                                                    & 0.000 & 0.000 & 0.005 &0.007 \\
 & &  $\phi$                                                              &  -0.005 & 0.011 &0.106  & 0.368\\[0.5ex]
 300 & (0.316,0.179)  & $\beta_{2}$                          & 0.033 & 0.041 &0.294  &0.336  \\
 & &  $\beta_{0}$                                                   &  -0.076 & -0.117 & 0.364 & 0.468 \\
  & &  $\beta_{1}$                                                   & 0.053 & 0.114 &0.219  &0.421  \\
 & &  $\gamma_1$                                                    &  -0.006& -0.006 &0.028  &0.027  \\ 
 & &  $\gamma_2$                                                    &  0.000& 0.000 &0.010  &0.011 \\
 & &  $\phi$                                                              &  0.036 & 0.062 &0.275  & 0.410\\[0.5ex]
 400 & (0.316,0.179)  & $\beta_{2}$                          & 0.033 & 0.027 & 0.245 & 0.289 \\
 & &  $\beta_{0}$                                                   &  -0.039 & -0.071 &0.288  &0.377  \\
  & &  $\beta_{1}$                                                   & 0.028 & 0.066 & 0.144 &  0.268\\
 & &  $\gamma_1$                                                    & -0.006 & -0.004 & 0.024 &  0.025\\ 
 & &  $\gamma_2$                                                    &  0.001& 0.000 & 0.009 &0.011 \\
 & &  $\phi$                                                              &  0.025 & 0.025 &0.218  & 0.410\\ 
\hline
\end{tabular}
\label{table:phi0.5}
\end{table} 
\begin{table} [ht!]
\caption{Comparison of PNCG algorithm with EM algorithm in terms of bias and RMSE with true value of $\phi$ as 0.75}
\centering
\begin{tabular}{ l l l l l l l }\\ 
\hline                                 
$n$ & $(\gamma_1,\gamma_2)$ & Parameter &\multicolumn{2}{c}{Bias}  &\multicolumn{2}{c}{RMSE} \\  \cline{4-5}  \cline{6-7}
& & &PNCG & EM & PNCG & EM \\
\hline  

300 & (0.215,0.183)  & $\beta_{2}$                          &0.018  & 0.013 & 0.247 & 0.368 \\
 & &  $\beta_{0}$                                                   & -0.050  & -0.138 &0.263  & 0.486 \\
  & &  $\beta_{1}$                                                   & 0.027 & 0.077 & 0.153 &0.354  \\
 & &  $\gamma_1$                                                    & -0.004 & -0.001 &0.020  & 0.018 \\ 
 & &  $\gamma_2$                                                    & 0.001 & 0.002 &0.007  & 0.008 \\
 & &  $\phi$                                                              & 0.000  & -0.140 & 0.159 & 0.422\\[0.5ex]
 400 & (0.215,0.183)  & $\beta_{2}$                          & 0.002 & -0.018 &0.195  &0.318  \\
 & &  $\beta_{0}$                                                   &  -0.029 & -0.094 &0.209  &0.417  \\
  & &  $\beta_{1}$                                                   & 0.025 & 0.043 &0.115  &0.231  \\
 & &  $\gamma_1$                                                    & -0.002 & 0.000 &0.017  & 0.017 \\ 
 & &  $\gamma_2$                                                    &  0.001& 0.002 & 0.006 &0.007 \\
 & &  $\phi$                                                              &  -0.003 & -0.138 & 0.119 & 0.408\\[0.5ex]
 300 & (0.316,0.179)  & $\beta_{2}$                          & 0.029 & -0.002 &0.300  &0.352  \\
 & &  $\beta_{0}$                                                   &  -0.063 & -0.143&0.370  &0.458  \\
  & &  $\beta_{1}$                                                   & 0.039 & 0.056 &0.193  &0.290  \\
 & &  $\gamma_1$                                                    &  -0.006& -0.003 & 0.029 &0.027  \\ 
 & &  $\gamma_2$                                                    &  0.001& 0.003 & 0.010 &0.012 \\
 & &  $\phi$                                                              &  0.006 &  -0.168& 0.274 & 0.445 \\[0.5ex]
 400 & (0.316,0.179)  & $\beta_{2}$                          &  0.034& 0.003 & 0.269 & 0.312  \\
 & &  $\beta_{0}$                                                   &  -0.045 & 0.077 & 0.315 & 0.385 \\
  & &  $\beta_{1}$                                                   & 0.028 &  0.055 & 0.166  &0.265  \\
 & &  $\gamma_1$                                                    & -0.004 & -0.004 &0.025  &0.025  \\ 
 & &  $\gamma_2$                                                    &  0.001&  0.001& 0.010 & 0.010 \\
 & &  $\phi$                                                              &  0.022 & -0.101 &0.256  & 0.418 \\ 
\hline
\end{tabular}
\label{table:phi0.75}
\end{table}

Next, we compare the performance of the PNCG algorithm, in terms of bias and RMSE, with the following three optimization routines readily available in R software: (i) "optim" with the "CG" (conjugate gradient) method based on the work by Fletcher and Reeves (1964), (ii) "nlm" that uses a Newton-type algorithm (Schnabel et al., 1985), and (iii) "Rcgmin" which is a conjugate gradient optimization that uses the Dai-Yuan update (Dai and Yuan, 2001). As in the case of PNCG, for the aforementioned optimization routines, maximization was also done on the observed log-likelihood function. Furthermore, we used the same technique of finding the initial values as we did for the PNCG method.

The implementation of any non-linear conjugate gradient scheme involves a line search algorithm and determining the search directions for descent of the function under consideration. Convergence of the non-linear conjugate gradient depends on the accuracy of the line search algorithm and the corresponding update of the search directions. Since exact line search algorithms are difficult to realize in practice, inexact line search algorithms are used based on the Wolfe conditions (see Wolfe, 1969). With the output of the step size obtained from a line search algorithm, the search directions are updated to obtain descent of the function. An accurate line search algorithm ensures that these search directions yield descent. Thus, it is important to devise a non-linear conjugate gradient algorithm that would provide descent directions in combination with inexact line search algorithms. It has been shown that the conjugate gradient method proposed by Fletcher and Reeves (1964) may not yield descent directions with inexact line search algorithms that satisfy the Wolfe conditions even for strongly convex problems. Also the scheme by Dai and Yuan (2001) requires a stronger version of the Wolfe conditions to prove convergence of the algorithm (Hager and Zhang, 2005; Hager and Zhang, 2006). However,
for general inexact line search algorithms, the Dai-Yuan method cannot be shown to be convergent. Our PNCG scheme with the update parameter $\xi_k$, originally proposed by Hager and Zhang (2005), counters these two major issues as it is convergent for any inexact line search algorithm and is more fast and accurate in comparison to the other two conjugate gradient algorithms. Thus the proposed PNCG is efficient and more robust with applications to a large class of non-linear, non-convex functions.

\subsection{Comparison with available optimization routines in R software}

For the purpose of this comparison, we consider two different settings - Setting 1: $\gamma_1=0.316,\gamma_2=0.179,\phi=0.75$ and Setting 2: $\gamma_1=0.316,\gamma_2=0.179,\phi=0.50$. In Table \ref{table:O1}, we present the comparison in terms of bias and RMSE. Note that the "optim", "nlm", and "Rcgmin" methods also allowed simultaneous estimation of all model parameters as we have seen in our proposed PNCG algorithm. However, our proposed PNCG algorithm outperforms other optimization techniques both in terms of bias and RMSE. Note, in particular, the reduction in bias and RMSE that the PNCG results with respect to the DNB shape parameter $\phi$. It is important to mention here that for both "optim" (method=CG) and "Rcgmin", the two conjugate gradient optimization techniques available in R, the number of iterations reached the maximum (set at 500) for every Monte Carlo run. As such, the reported bias and RMSE for these two methods actually correspond to the estimates obtained at the last iteration (500th) step. This is unlike our proposed PNCG algorithm, where the convergence took place in less than 85 iterations (on an average) for any considered parameter setting. The "nlm" routine, on the other hand, took less than 15 iterations (on an average) to converge. Even though the "nlm" converges faster than PNCG, the final estimates obtained from PNCG is more accurate and more precise (as already seen from Table \ref{table:O1}). 

\begin{table} [ht!]
\caption{Comparison of PNCG algorithm with other optimization algorithms}
\centering
\begin{tabular}{ l l l l l l l || l l l l}\\ 
\hline                                 
$n$ & Setting & Parameter &\multicolumn{4}{c}{Bias}  &\multicolumn{4}{c}{RMSE} \\  \cline{4-7}  \cline{8-11}
& & &PNCG & optim &nlm &Rcgmin & PNCG & optim &nlm &Rcgmin \\
\hline  

300 &1 & $\beta_{2}$                                             &0.032  & 0.056 & 0.047 & 0.203        &0.303 & 0.357&0.303 & 0.525 \\
 & &  $\beta_{0}$                                                    & -0.094  & -0.118 &-0.066  & -0.016   &0.388 &0.462 &0.384 &1.050 \\
  & &  $\beta_{1}$                                                   & 0.059 & 0.084 & 0.064 &0.277         &0.221 &0.266 &0.316 &1.140 \\
 & &  $\gamma_1$                                                 & -0.006 & -0.007 &-0.005  & -0.017    &0.030 &0.031 &0.031 & 0.039\\ 
 & &  $\gamma_2$                                                 & 0.001 & 0.001 &0.001  & -0.003        &0.012 &0.012 &0.011 & 0.014\\
 & &  $\phi$                                                           & 0.024  & -0.069 & 0.087 & 0.551       &0.324 &0.505 &0.496 & 1.155\\[1ex]
 
 400 &1 & $\beta_{2}$                                            &0.014  & 0.046 & 0.029 & 0.152       &0.271 &0.320 &0.243 &0.400 \\
 & &  $\beta_{0}$                                                    & -0.031  & -0.061 &-0.005  & 0.018   &0.305 &0.398 &0.302 &0.594 \\
  & &  $\beta_{1}$                                                   & 0.031 & 0.063 & 0.028 &0.135         &0.157 &0.227 &0.157 &0.407 \\
 & &  $\gamma_1$                                                  & -0.001 & -0.002 &-0.002  & -0.010   &0.024 &0.026 &0.025 &0.031 \\ 
 & &  $\gamma_2$                                                 & 0.000 & -0.001 &0.000  & -0.004      &0.009 &0.010 &0.010 &0.012 \\
 & &  $\phi$                                                             & -0.003  & -0.061 & 0.067 & 0.425     &0.249 &0.475 &0.407 & 0.984\\[1ex]
 
 300 &2 & $\beta_{2}$                                             &0.033  & 0.086 & 0.048 & 0.208        &0.285 &0.359 &0.265 &0.469 \\
 & &  $\beta_{0}$                                                    & -0.042  & -0.061 &-0.032  & -0.027    &0.346 &0.448 &0.341 &0.690 \\
  & &  $\beta_{1}$                                                   & 0.045 & 0.105 & 0.055 &0.284          &0.181 &0.274 &0.282 &1.263 \\
 & &  $\gamma_1$                                                 & -0.007 & -0.012 &-0.007  & -0.020     &0.030 &0.033 &0.032 &0.039 \\ 
 & &  $\gamma_2$                                                 & 0.001 & -0.001 &0.000  & -0.005       &0.011 &0.011 &0.010 &0.013 \\
 & &  $\phi$                                                           & 0.030  & 0.212 & 0.077 & 0.595         &0.238 &0.517 &0.386 & 0.981\\[1ex]
 
 400 &2 & $\beta_{2}$                                            &0.041  & 0.085 & 0.035 & 0.167        &0.263 &0.327 &0.228 & 0.413\\
 & &  $\beta_{0}$                                                    & -0.045  & -0.053 &-0.033  & -0.018   &0.316 &0.419 &0.363 & 0.716\\
  & &  $\beta_{1}$                                                   & 0.029 & 0.067 & 0.045 &0.194          &0.167 &0.245 &0.353 & 0.670\\
 & &  $\gamma_1$                                                  & -0.004 & -0.008 &-0.005  & -0.014    &0.024 &0.026 &0.025 & 0.031\\ 
 & &  $\gamma_2$                                                 & 0.000 & -0.001 &0.000 & -0.004        &0.010 &0.011 &0.010 & 0.012\\
 & &  $\phi$                                                             & 0.015  & 0.148 & 0.050 & 0.437       &0.216 &0.459 &0.290 & 0.850\\[1ex]
 
\hline
\end{tabular}
\label{table:O1}
\end{table}

\section{Real data analysis}

To illustrate the proposed PNCG algorithm, we analyzed a well-known melanoma data that is available in the ``timereg'' package of R software. The dataset contains 205 patients who were observed after operation for removal of malignant melanoma in the period 1962-1977 and then followed until 1977. The observed time is recorded in years and refers to the time until the patient's death or the censoring time. It has a mean of 5.9 years and a standard deviation of 3.1 years. Patients who survived the end of the study and patients who died due to some other causes were considered as censored observations. The percentage of censored observations is 72\%. For our analysis, we selected the ulceration status (presence of ulcer for 90 patients and absence of ulcer for 115 patients) and tumor thickness (measured in mm with a mean of 2.92 and a standard deviation of 2.96) as prognostic factors. As done in Pal and Balakrishnan (2016), we linked the parameter $p$ to tumor thickness and the parameter $\eta$ to ulceration status to retain identifiability of model parameters; see also Rodrigues et al. (2011). 

We employed the proposed PNCG algorithm and the results are presented in Table \ref{table:MLE}. For comparison purpose, we also present the EM, optim, nlm, and Rcgmin estimates together with their standard errors. For the EM algorithm, the values are taken from Table 3 of Pal and Balakrishnan (2016). Note that in the EM algorithm, the M-step was carried out using a one-step Newton Raphson method and the estimates of $\gamma_1$ and $\phi$ were obtained by employing a two-way profile likelihood approach within the EM algorithm. Furthermore, the standard errors corresponding to the EM estimates were obtained by inverting the observed information matrix under the assumption that $\gamma_1$ and $\phi$ were fixed. For the optim and nlm routines, the hessian matrix is readily available as an output which we then inverted to calculate the standard errors. For the Rcgmin routine, however, we computed the hessian matrix numerically using the package "numDeriv" and using the "Richardson" method. For the destructive model that we considered, we also tried to compute the standard errors by inverting the observed information matrix (with respect to all model parameters), however, we noticed that the second-order derivatives of the observed log-likelihood function were highly unstable, specifically with respect to the parameter $\phi.$ As such, corresponding to the PNCG estimates, we calculated the standard errors using the non-parametric bootstrap method. It is clear that the bootstrap standard errors are much lower compared to the EM, optim, nlm, and Rcgmin standard errors. Note also that the PNCG, EM, optim, and nlm estimates are relatively close to each other.
\begin{table} [ht!]
\caption{MLEs and standard errors of the parameters for the melanoma data.}
\begin{tabular*}{\textwidth}{@{\extracolsep{\fill}} c c c c c c || c c c c c}\\
\hline                                 
Parameter                          &\multicolumn{5}{c}{Estimate} &\multicolumn{5}{c}{Standard Error} \\  \cline{2-6}\cline{7-11}
& PNCG & EM & optim &nlm &Rcgmin & PNCG & EM & optim &nlm &Rcgmin    \\
\hline                                 
$\beta_{1,\text {intercept}}$                            &  -5.841 &-5.882  &  -5.787   & -5.880 & -6.434 &                   0.573 & 1.183    & 2.099  & 2.155 & 2.473  \\
$\beta_{1,\text{thickness}}$                             &  1.183 &  1.197 &   1.191    & 1.190  & 1.284 &                 0.383 & 0.354     &  0.444  & 0.451 &  0.499 \\
$\beta_{2,\text{ulc:present}}$                          &  5.434 &5.490  &   5.536  & 5.490 & 6.432 &                     0.523 & 1.001    &  2.651 & 2.675 &   3.186  \\
$\beta_{2,\text{ulc:absent}}$                            &  3.533 &3.484 &   3.523    & 3.480 &  4.256  &                    0.479 &  1.154  &  2.247  & 2.266 &  2.672  \\
$\gamma_1$                                                   &   0.314 & 0.300 & 0.308    & 0.311 & 0.284 &                      0.029 &      -     &  0.080 & 0.081 &  0.081         \\
$\gamma_2$                                                   &  0.122 &0.127 &   0.122   & 0.123 & 0.119 &                       0.020 & 0.021  & 0.026  & 0.026 &   0.025    \\
$\phi$                                                             & 6.654 & 6.600&    7.146   & 6.600 & 8.340 &                        0.604 &      -     &  3.525 & 3.205 &   4.274             \\
\hline
\end{tabular*}
\label{table:MLE}
\end{table}
 
\section{Concluding remarks}

In this paper, we have proposed an algorithm based on a projected non-linear conjugate gradient method with an efficient line search technique (Hager and Zhang, 2005) for the maximum likelihood estimation of the destructive cure rate model parameters. For simulation studies, we assumed the initial number of competing risks to follow a negative binomial distribution and showed that our proposed method results in estimates that have smaller bias and smaller RMSE when compared to the estimates obtained from the EM algorithm as well as other existing Newton-type and conjugate gradient type methods. The reduction in the RMSE is more pronounced for the parameters that are associated with the cure rate, which holds the key for accurate inference on overall population survival. Our proposed algorithm, being a general algorithm, can be used for destructive cure rate models with any competing risks distribution. Furthermore, on setting the activation probability $p$ to zero, the destructive cure rate model reduces to the regular cure rate model (Berkson and Gage, 1952; Chen et al., 1999) and as such our algorithm can be used to study these regular cure models as well. We believe that we are the first one to use the PNCG algorithm with efficient line search technique in the context of cure rate models and we hope that researchers will use our algorithm for more complicated cure rate models, for instance, those considered in Gallardo et al. (2016) and Pal et al. (2018), among others. A possible extension will be to consider interval censoring, a more general form of censoring, and develop the likelihood inference for destructive cure models based on the PNCG algorithm. We are currently working on this and hope to report the findings in a future paper.
  
\section*{Conflict of interest}

On behalf of all authors, the corresponding author states that there is no conflict of interest.

\section*{Acknowledgment}

Both authors express their thanks to the National Institutes of Health Grant Number 1R21CA242933-01 to support this research. 

\section*{References}

\begin{description}

\item  Adesokan, B., Knudsen, K., Krishnan, V. P. and Roy, S. (2018) A fully non-linear optimization approach to acousto-electric tomography. {\it Inverse Problems} {\bf 34}, 104004.

\item Annunziato, M. and Borz\`i, A. (2013) A Fokker–Planck control framework for multidimensional stochastic process. {\it Journal of Computational and Applied Mathematics} {\bf 237}, 487-507.

\item Balakrishnan, N., Barui, S., and Milienos, F. (2017). Proportional hazards under Conway–Maxwell-Poisson cure rate model and associated inference. {\it Statistical Methods in Medical Research} {\bf 26}, 2055$–$2077.

\item Balakrishnan, N. and Feng, T. (2018). Proportional odds under Conway-Maxwell-Poisson cure rate model and associated likelihood inference. {\it Statistics, Optimization and Information Computing} {\bf 6}, 305$-$334.

\item Balakrishnan, N., Koutras, M. V., Milienos, F. S., and Pal, S. (2016). Piecewise linear approximations for cure rate models and associated inferential issues. {\it Methodology and Computing in Applied Probability} {\bf 18}, 937$-$966.

\item Balakrishnan, N. and Pal, S. (2013). Lognormal lifetimes and likelihood-based inference for flexible cure rate models based on COM-Poisson family. {\it Computational Statistics \& Data Analysis} {\bf 67}, 41$-$67.

\item Balakrishnan, N. and Pal, S. (2015). An EM algorithm for the estimation of parameters of a flexible cure rate model with generalized gamma lifetime and model discrimination using likelihood- and
information-based methods. {\it Computational Statistics} {\bf 30}, 151$-$189.

\item Balakrishnan, N. and Pal, S. (2016). Expectation maximization-based likelihood inference for flexible cure rate models with Weibull lifetimes. {\it Statistical Methods in Medical Research} {\bf 25}, 1535$-$1563.

\item Berkson, J. and Gage, R. P. (1952). Survival curve for cancer patients following treatment. {\it Journal of the American Statistical Association} $\boldsymbol{47}$, 501$-$515.

\item Boag, J. W. (1949). Maximum likelihood estimates of the proportion of patients cured by cancer therapy. {\it Journal of the Royal Statistical Society: Series B} $\boldsymbol{11},$ 15$-$53.


\item Cancho, V. G., Ortega, E. M. M., Barriga, G. D. C. and Hashimoto, E. M. (2011). The Conway-Maxwell-Poisson-generalized gamma regression model with long-term survivors. {\it Journal of Statistical Computation and Simulation} {\bf 81}, 1461$-$1481.

\item Cancho, V. G., Bandyopadhyay, D. and Louzada, F. (2013). The destructive negative binomial cure rate model with a latent activation scheme. {\it Statistical Methodology} {\bf 13}, 48$-$68.

\item Chen, M. -H., Ibrahim, J. G. and Sinha, D. (1999). A new Bayesian model for survival data with a surviving fraction. {\it Journal of the American Statistical Association} $\boldsymbol{94}$, 909$-$919.

\item Cooner, F., Banerjee, S., Carlin, B. P. and Sinha, D. (2007). Flexible cure rate modeling under latent activation schemes. {\it Journal of the American Statistical Association} $\bf{102}$, 560$-$572.

\item Dai, Y. H. and Yuan, Y. (2001). An efficient hybrid conjugate gradient method for unconstrained optimization. {\it Annals of Operations Research} {\bf 103}, 33$-$47.

\item de Freitas, L. A. and Rodrigues, J. (2013). Standard exponential cure rate model with informative censoring. {\it Communications in Statistics - Simulation and Computation} {\bf 42}, 8$–$23.


\item Dunn, P. K. and Smyth, G. K. (1996). Randomized quantile residuals. {\it Journal of Computational and Graphical Statistics} {\bf 5}, 236$–$244.

\item Farewell, V. T. (1986). Mixture models in survival analysis: are they worth the risk? {\it Canadian Journal of Statistics} {\bf 14}, 257$-$262.

\item Fletcher, R. and Reeves, C. M. (1964). Function minimization by conjugate gradients. {\it Computer Journal} 7, 148$-$154.

\item Gallardo, D. I., Bolfarine, H. and Pedroso-de Lima, A. C. (2016). An EM algorithm for estimating the destructive weighted Poisson cure rate model. {\it Journal of Statistical Computation and Simulation} {\bf 86}, 1497$-$1515.

\item Hager, W. W. and Zhang, H. (2005). A new conjugate gradient method with guaranteed descent and an efficient line search. {\it SIAM Journal on Optimization} {\bf 16}, 170-192.

\item Hager, W. W. and Zhang, H. (2006).  A survey of nonlinear conjugate gradient methods. {\it Pacific Journal of Optimization} {\bf 2}, 35$-$58.


\item Kannan, N., Kundu, D., Nair, P. and Tripathi, R. C. (2010). The generalized exponential cure rate model with covariates. {\it Journal of Applied Statistics} {\bf 37}, 1625$-$1636.


\item Kuk, A. Y. C. and Chen, C. H. (1992). A mixture model combining logistic regression with proportional hazards regression. {\it Biometrika} {\bf 79}, 531$-$541.

\item Li, C. S. and Taylor, J. M. G. (2002). A semi-parametric accelerated failure time cure model. {\it Statistics in Medicine} {\bf 21}, 3235$–$3247.

\item Li, C. S., Taylor, J. M. G. and Sy, J. P. (2001). Identifiability of cure models. {\it Statistics and Probability Letters} {\bf 54}, 389$-$395.

\item Maller, R. A. and Zhou, X. (1996). {\it Survival Analysis with Long-Term Survivors}. New York: John Wiley \& Sons.

\item McLachlan, G. J. and Krishnan, T. (2008). {\it The EM Algorithm and Extensions}, Second edition. Hoboken, New Jersey: John Wiley \& Sons.

\item Neittaanmaki, P. and Tiba, D. (1994). Optimal Control of Nonlinear Parabolic Systems: Theory, Algorithms and Applications. {\it Pure and Applied Mathematics. CRC Press, London}.

\item Nocedal, J. and Wright, S. J. (1999). {\it Numerical Optimization}. New York: Springer.


\item Pal, S. and Balakrishnan, N. (2016). Destructive negative binomial cure rate model and EM-based likelihood inference under Weibull lifetime. {\it Statistics and Probability Letters} {\bf 116}, 9$-$20.

\item Pal, S. and Balakrishnan, N. (2017). Likelihood inference for the destructive exponentially weighted Poisson cure rate model with Weibull lifetime and an application to melanoma data. {\it Computational Statistics} {\bf 32}, 429$-$449.

\item Pal, S., Majakwara, J. and Balakrishnan, N. (2018). An EM algorithm for the destructive COM-Poisson regression cure rate model. {\it Metrika} {\bf 81}, 143$-$171.


\item Rodrigues, J., de Castro, M., Balakrishnan, N. and Cancho, V. G. (2011). Destructive weighted Poisson cure rate models. {\it Lifetime Data Analysis} $\boldsymbol{17},$ 333$-$346.

\item Rodrigues, J., de Castro, M., Cancho, V. G. and Balakrishnan, N. (2009). COM-Poisson cure rate survival models and an application to a cutaneous melanoma data. {\it Journal of Statistical Planning and Inference} $\boldsymbol{139},$ 3605$-$3611.

\item Roy, S., Annunziato, M. and Borz\`i, A. (2016). A Fokker–Planck feedback control-constrained approach for modelling crowd motion. {\it Journal of Computational and Theoretical Transport} {\bf 45}, 442-458.

\item Roy, S., Annunziato, M., Borz\`i, A. and Klingenberg, C. (2018). A Fokker–Planck approach to control collective motion. {\it Computational Optimization and Applications} {\bf 69}, 423-459.

\item Roy, S., Borz\`i, A. and Habbal, A. (2017). Pedestrian motion constrained by FP-constrained Nash games. {\it Royal Society Open Science}, {\bf 4}, 170648.


\item Schnabel, R. B., Koontz, J. E. and Weiss, B. E. (1985). A modular system of algorithms for unconstrained minimization. {\it ACM Transactions on Mathematical Software} {\bf 11}, 419$-$440.

\item Sy, J. P. and Taylor, J. M. G. (2000). Estimation in a Cox proportional hazards cure model. {\it Biometrics} $\boldsymbol{56},$ 227$-$236.

\item Wolfe, P. (1969). Convergence conditions for ascent methods. {\it SIAM Review} {\bf 11}, 226$-$235.

\item Yang, G. and Chen, C. (1991). A stochastic two-stage carcinogenesis model: a new approach to computing the probability of observing tumor in animal bioassays. {\it Mathematical Biosciences} {\bf 104}, 247$-$258.



\end{description}

\end{document}